\documentclass[12pt]{article}

\usepackage{amssymb,amsmath,amsthm, amsfonts}

\newcommand{\f}{\mathcal{F}}
\newcommand{\e}{\mathcal{E}}
\newcommand{\lft}{L_{+}^2(\mathcal{F}_T)}
\newcommand{\lf}{L^{2}(\mathcal{F}_T)}
\newcommand{\lfst}{L^{2}(\mathcal{F}_T;S_T)}

\def\sqr#1#2{{\vcenter{\vbox{\hrule height.#2pt
              \hbox{\vrule width.#2pt height#1pt \kern#1pt \vrule width.#2pt}
              \hrule height.#2pt}}}}
\def\signed #1{{\unskip\nobreak\hfil\penalty50
              \hskip2em\hbox{}\nobreak\hfil#1
              \parfillskip=0pt \finalhyphendemerits=0 \par}}
\def\endpf{\signed {$\sqr69$}}
\def\3n{\negthinspace \negthinspace \negthinspace }
\def\2n{\negthinspace \negthinspace }
\def\1n{\negthinspace }

\def\={\buildrel \triangle \over =}

\def\a{\alpha}

\def\d{\delta}
\def\O{\Omega}

\def\cA{{\cal A}}

\def\cE{{\cal E}}
\def\cF{{\cal F}}

\def\cL{{\cal L}}

\def\no{\noindent}

\def\ms{\medskip}
\def\bs{\bigskip}
\def\q{\quad}
\def\qq{\qquad}

\def\exp{\mathop{\rm exp}}
\def\sup{\mathop{\rm sup}}

\def\ae{\hbox{\rm a.e.{ }}}
\def\as{\hbox{\rm a.s.{ }}}

\def\|{\Big |}
\def\({\Big (}
\def\){\Big )}
\def\[{\Big[}
\def\]{\Big]}
\def\be{\begin{equation}}
\def\bel{\begin{equation}\label}
\def\ee{\end{equation}}
\def\bt{\begin{theorem}}
\def\bcd{\begin{condition}}
\def\ecd{\end{condition}}
\def\et{\end{theorem}}
\def\bc{\begin{corollary}}
\def\ec{\end{corollary}}
\def\bde{\begin{definition}}
\def\ede{\end{definition}}
\def\bl{\begin{lemma}}
\def\el{\end{lemma}}
\def\bp{\begin{proposition}}
\def\ep{\end{proposition}}
\def\br{\begin{remark}}
\def\er{\end{remark}}
\def\ba{\begin{array}}
\def\ea{\end{array}}
\def\ed{\end{document}}

\def\square#1{\vbox{\hrule\hbox{\vrule height#1%
     \kern#1\vrule}\hrule}}
\def\rectangle#1#2{\vbox{\hrule\hbox{\vrule height#1%
     \kern#2\vrule}\hrule}}

\font\tenbb=msbm10 \font\sevenbb=msbm7 \font\fivebb=msbm5

\newfam\bbfam
\scriptscriptfont\bbfam=\fivebb \textfont\bbfam=\tenbb
\scriptfont\bbfam=\sevenbb

\newtheorem{lemma}{Lemma}[section]
\newtheorem{remark}{Remark}[section]

\newtheorem{theorem}{Theorem}[section]
\newtheorem{corollary}{Corollary}[section]

\newtheorem{definition}{Definition}[section]
\newtheorem{proposition}{Proposition}[section]
\newtheorem{condition}{Condition}[section]

\makeatletter
   
   \@addtoreset{equation}{section}
\makeatother

\begin{document}

    \title{{\bf  Filtration-Consistent Dynamic Operator with a Floor and  Associated Reflected Backward Stochastic Differential
    Equations}\thanks{This
work is partially supported by the NSFC under grants 10325101
(distinguished youth foundation) and 101310310 (key project), and
the Science Foundation of Chinese Ministry of Education under
grant 20030246004.}}

\author{Xiaobo Bao\thanks{Institute of Mathematics, School of Mathematical Sciences, Fudan University, Shanghai
200433, China, \& Key Laboratory of Mathematics  for Nonlinear
Sciences (Fudan University), Ministry of Education.}  \quad and
\quad Shanjian Tang\thanks{Department of Finance and Control
Sciences, School of Mathematical Sciences, Fudan University,
Shanghai 200433, China, \& Key Laboratory of Mathematics  for
Nonlinear Sciences (Fudan University), Ministry of Education.
{\small\it E-mail:} {\small\tt sjtang@fudan.edu.cn}.\ms} }

\date{}
\maketitle

\abstract{This paper introduces the notion of an $\{\cF_t, 0\le
t\le T\}$-consistent  dynamic operator with a floor, by suitably
formulating four axioms. It is shown that  an $\{\cF_t, 0\le t\le
T\}$-consistent   dynamic operator $\{\cE_{s,t}, 0\le s\le t\le
T\}$ with a continuous upper-bounded floor $S$ is necessarily
represented by the solution of a backward stochastic differential
equation (BSDE) reflected upwards on the floor, if it is
$\cE^\mu$-super-dominated for some $\mu>0$ and if it has the
non-increasing and floor-above-invariant property of forward
translation. We make full use of the two assumptions to extend the
underlying $\{\cF_t, 0\le t\le T\}$-consistent dynamic operator
from the subset $L^2(\cF_T;S_T)$ of floor-dominating
square-integrable random variables to the whole space $L^2(\cF_T)$
of square-integrable random variables. The extended dynamic
operators are shown to be identified to an $\{\cF_t, 0\le t\le
T\}$-consistent expectation $\widetilde\cE$. The generator $g$ of
its BSDE representation given by Coquet, Hu, Memin, and
Peng~\cite[\it Probability Theory and related Fields, 123 (2002),
1--27]{CHMP} turns out to be that of the desired RBSDE. More
elegantly, the process $\{\cE_{t,T}[X], 0\le t\le T\}$ turns out
to be an $\widetilde \cE$-supermartingale for each $X\in
L^2(\cF_T;S_T)$. This key observation allows us to apply a
nonlinear Doob-Meyer's decomposition theorem to construct the
increasing process of $\{\cE_{t,T}[X], 0\le t\le T\}$ as a
solution to some BSDE reflected upwards on the floor. }

\bs

\no{\bf 2000 Mathematics Subject Classification}.  Primary 60H10;
Secondary 60H30.

\bs

\no{\bf Key Words}.  reflected BSDE, dynamic nonlinear operator,
filtration-consistent nonlinear expectation, nonlinear Doob-Meyer
decomposition.

\section{Introduction.}

  Let $\{B_t,0 \leq t \leq T\}$ be a $d$-dimensional standard
Brownian motion defined on a probability space
$(\Omega,\mathcal{F},P)$. Let $\{\mathcal{F}_t,0 \leq t \leq T\}$
be the natural filtration of $\{B_t, t\in [0,T]\}$, augmented by
all $P$-null sets of $\mathcal{F}$.
  Set
  $$\ba{rcl} L^2(\mathcal{F}_T)&:=&\{ \xi: \xi \hbox{ \rm is an $\mathcal{F}_T$-measurable random variable  s.t. } E|\xi|^2<
  +\infty\},\\
L^{2}(\mathcal{F}_T;S_T)&:=&\{ \xi \in L^2(\mathcal{F}_T): \xi
\geq S_T\},\quad  L_{+}^2(\mathcal{F}_T):=L^2(\mathcal{F}_T;0),\\
 \cL_{\cF}^2(0,T;{\mathbb{R}}^m)&:=&\{\phi:  \hbox{ \rm $\phi$ is $\mathbb{R}^m$-valued and  $\{\cF_t, 0\le t\le T\}$-adapted
 s.t. }  E\int_0^T|\varphi_t|^2dt<\infty\},\\
  \mathbb{H}^2&:=& \Big\{\varphi\in \cL_\cF^2(0,T;\mathbb{R}): \varphi \hbox{ \rm is predictable
  }\Big\},\\
  \mathbb{S}^2&:=& \Big\{\varphi\in \mathbb{H}^2: \varphi \hbox{ \rm is
  right-continuous  s.t. } E\sup_{0 \leq t \leq T}|\varphi_t|^2<
  +\infty\Big\}.\ea$$

Let be given three objects:  a terminal value $\xi$, a random
field $g:\Omega \times[0,T] \times \mathbb{R}\times
  \mathbb{R}^d\longrightarrow \mathbb{R}$, and a continuous progressively measurable real-valued
  random process $\{S_t,0\leq t \leq T\}$. Assume that

(C1) $\xi \in L^2(\cF_T)$.

(C2)  $g(\cdot,y,z)\in
  \mathbb{H}^2 \hbox{ \rm for } (y,z)\in \mathbb{R} \times
  \mathbb{R}^d$.

  (C3) $|g(t,y,z)-g(t,y^{'},z^{'})|\leq \a (|y-y^{'}|+|z-z^{'}|)$ a.s. with $y,y^{'}\in \mathbb{R}$ and $z,z^{'}\in \mathbb{R}^d
  $ for some positive constant $\a$. And

    (C4) $\displaystyle S^+\in \mathbb{S}^2$ and  $S_T\leq \xi$ a.s..

Consider the following reflected backward stochastic differential
equation
    (RBSDE):
 \be \left\{\ba{c}\displaystyle Y_t=\xi + \int_t^T g(s,Y_s,Z_s)ds+ K_T-K_t-\int_t^T\langle Z_s,d
    B_s\rangle ,\qq 0\leq t \leq T;\\
\displaystyle Y_t\geq S_t, a.s. 0\leq t \leq T; \q K_0=0 \hbox{
\rm and }
    \int_0^T (Y_t-S_t)d K_t=0.\ea\right.\label{RBSDE0}
\ee Its  solution  is a triple $\{(Y_t,Z_t,K_t),0\leq t \leq
    T\}$ of $\{\mathcal{F}_t, 0\le t\le T\}$-progressively measurable processes
    taking values in $\mathbb{R}\times \mathbb{R}^d\times
    \mathbb{R}_{+}$ such that (i) $Z\in \mathbb{H}^2$,
and (ii) $\{K_t, 0\le t\le T\}$ is continuous and increasing.
There exists unique solution $\{(Y_t,Z_t,K_t),0\leq t \leq T\}$ of
RBSDE~(\ref{RBSDE0}) if conditions (C1)-(C4) are satisfied. For
further details, see El Karoui,
    Kapoudjian, Pardoux, Peng, and Quenez~\cite{ELKa}.

Define  $$\e_{t,T}^{r;g,S}[Y]:=Y_t,\q  \forall\, Y \in\lfst$$
  and
 $$\e_{t,T}^{r,g}[Y]:=\e_{t,T}^{r;g,0}[Y],\q \forall Y\, \in\lft.$$
 Here, the superscript $r$ indicates that the underlying operators
 are generated by a reflected BSDE,
 and the superscripts  $(g,S)$ specify the generator and the obstacle of the underlying RBSDE.
 It was noted by El Karoui and Quenez~\cite{ELKaQuenez} that the
 operator $\e_{t,T}^{r;g,S}$ satisfies the following properties:

($\cA$1) {\bf Monotonicity.} $\mathcal{E}_{s,t}[Y]\geq
\mathcal{E}_{s,t}[Y^{'}]$ if $Y,Y'\in L^2(\cF_t)$ and $Y\geq
Y^{'}$.

($\cA$2) {\bf Time consistency.}
$\mathcal{E}_{r,s}[\mathcal{E}_{s,t}[Y]]=\mathcal{E}_{r,t}[Y]$ if
$r\leq s\leq t\leq T$ and $Y\in L^2(\cF_t)$.
\\  Therefore it is a nonlinear pricing system of
 square-integrable American contingent claims. A financial
 theoretically interesting problem is the converse one: Is a
 dynamic operator with similar properties like ($\cA$1) and ($\cA$2) necessarily represented by a
 RBSDE? That is, we are concerned with axiomatic characteristics
 for a dynamic operator to be represented by a RBSDE.

This paper is concerned with the converse problem for
RBSDE~(\ref{RBSDE0}) with a general floor $S>-\infty$. More
precisely, we introduce a class of dynamic operators with floors,
and then relate them to BSDEs reflected upwards on the floors. In
this way, we characterize on one hand the solutions of RBSDEs with
axioms, and on the other hand, we give a representation for the
dynamic operators in terms of RBSDEs~(\ref{RBSDE0}).

Throughout the paper, we shall make  the following three
assumptions. The first one concerns some continuity of the
underlying dynamic operators:

 (H1)  {\bf
$\mathcal{E}^\mu$-super-domination} \it There is  some $\mu>0$
such that
\begin{equation}
\mathcal{E}_{t,T}[X+Y]-\mathcal{E}_{t,T}[X]\leq
\mathcal{E}_{t,T}^{\mu}[Y],\  \as
\end{equation}
for $ t\in [0,T], X\in L^{2}(\mathcal{F}_T;S_T) \hbox{ \rm and }
Y\in L_{+}^2(\mathcal{F}_T).$
 Here,  $\mathcal{E}^{\mu}_{t,T}[Y]:=y_t$ with  $\{y_t,0\leq t\leq T\}$ being the first component
 of the solution of the following BSDE:
  \be \label{bsde0} y_t=Y + \int_t^T \mu |z_s|ds-\int_t^T\langle z_s,d
    B_s\rangle ,\q 0\leq t \leq T.
 \ee\rm
The following assumption is its weaker version:

 (H1)'  {\bf
$\mathcal{E}^\mu$-domination} \it There is  some $\mu>0$ such that
for $X\in L^{2}(\mathcal{F}_T;S_T) \hbox{ \rm and } Y\in
L_{+}^2(\mathcal{F}_T)$, we have
\begin{equation}
\mathcal{E}_{0,T}[X+Y]-\mathcal{E}_{0,T}[X]\leq
\mathcal{E}_{0,T}^{\mu}[Y].
\end{equation}\rm

Note that the assumption (H1)' is much weaker than that of
$\mathcal{E}^\mu$-domination used by Coquet, Hu, Memin and
Peng~\cite[Definition 4.1, page 9]{CHMP}, in that $Y$ is here
restricted within $L^2_+(\cF_T)$ instead of taking values in the
whole space $L^2(\cF_T)$ like the latter. This difference will
complicate our subsequent arguments.

The second assumption still concerns the operators:

(H2) { \bf The non-increasing and floor-above-invariant property
of forward translation } $\mathcal{E}_{t,T}[X+Y]\leq
\mathcal{E}_{t,T}[X]+Y$, \as for $X \in L^{2}(\mathcal{F}_T;S_T)$
and $Y\in L_{+}^2(\mathcal{F}_t)$. Furthermore,  for $ X \in
L^{2}(\mathcal{F}_T;S_T)$ such that $\mathcal{E}_{t,T}[X]>S_t\ \as
$ for $ t\in [0,T]$, We have $\mathcal{E}_{t,T}[X+Y]=
\mathcal{E}_{t,T}[X]+Y$ for any $Y\in L_{+}^2(\mathcal{F}_t).$

The third assumption concerns the underlying floor $S$:

 (H3) \it The predictable process $S$ is continuous and $S^+\in \mathbb{S}^2$. Moreover, there is a positive constant   $C$  such
that $ \mbox{ess}\sup_{0\leq t\leq T}S_t\leq
 C $ a.s.. \rm

We now introduce the notion of an $\{\cF_t, 0\le t\le
T\}$-consistent  dynamic operator with floor $S$.

\bde \label{definef}
 A time parameterized  system of operators
\begin{equation}
\mathcal{E}_{s,t}[\cdot]:L^{2}(\mathcal{F}_t;S_t)
 \longrightarrow L^{2}(\mathcal{F}_s;S_s), \q 0\leq s\leq t\leq
 T
\end{equation}
is called an $\{\cF_t, 0\le t\le T\}$-consistent  dynamic operator
with floor $S$ if it satisfies : $\mathcal{E}_{s,t}[Y]$ is
continuous in $s\in[0, t]$ for $t\in [0,T]$ and $Y\in
L^{2}(\mathcal{F}_t;S_t)$, and furthermore, it satisfies the
following four axioms.

\rm (D1) \it { \bf Floor-above  strict monotonicity }
$\mathcal{E}_{s,t}[Y]\geq \mathcal{E}_{s,t}[Y^{'}]$ a.s. for $Y$
and $Y^{'}\in L^{2}(\mathcal{F}_t; S_t)$ such that $Y\geq Y^{'}$
a.s..  If $\e_{s,t}[Y]
>S_s\ \as \ \hbox{ \it for any } s\in [r,t]$ and some $r\in [0,t]$, then $Y'=Y$
a.s. if $Y'\geq Y$ a.s. and $\mathcal{E}_{r,t}[Y']=
\mathcal{E}_{r,t}[Y]$ a.s..

\rm (D2) \it $\mathcal{E}_{s,t}[Y]=Y$ a.s. for each $ Y\in
L^2(\mathcal{F}_s;C).$

\rm (D3) \it { \bf Time consistency }
$\mathcal{E}_{r,s}[\mathcal{E}_{s,t}[Y]]=\mathcal{E}_{r,t}[Y]$
a.s. for $Y\in L^2(\cF_t;C)$ if $r\leq s\leq t\leq T$.

\rm (D4) \it { \bf Zero-one law } For each $s\leq t$,
$\mathcal{E}_{s,t}[1_{A}Y+\tilde C]-\tilde
C=1_{A}(\mathcal{E}_{s,t}[Y+\tilde C]-\tilde C)$, a.s., $\forall
A\in \mathcal{F}_s$ for any constant $\tilde C$ dominating the
floor. \ede

The main result of the paper is stated as follows.

\bt \label{RBSDEConverse} Consider an $\{\cF_t, 0\le t\le
T\}$-consistent dynamic operator $\{\cE_{s,t}, 0\le s\le t\le T\}$
with floor $S$. Assume that $\{\cE_{s,t}, 0\le s\le t\le T\}$
satisfies (H1) and (H2), and  the obstacle $S$ satisfies (H3).
Then, there is a random field $g: \O\times [0, T]\times
\mathbb{R}^d\to \mathbb{R}$
 such that the following are satisfied:

 (i) $g(t,0)=0$ for $\ae t\in [0,T]$,

 (ii) $|g(t,z_1)-g(t,z_2)|\le \mu |z_1-z_2|$, and

 (iii) $\cE_{s,t}(Y)=\cE^{r;g,S}_{s,t}(Y)$ for  any $Y\in L^2(\cF_t;S_t)$ and $s\in [t,T]$ with $t\in [0,T]$.

For the particular case of the zero floor ($S=0$), the last
assertion is still true if the assumption (H1) is replaced by the
weaker one (H1)'. \et

It is worth noting that the extremal case of the negative infinite
floor (that is, $S=-\infty$) for RBSDE~(\ref{RBSDE0}) is the
following backward stochastic differential equation (BSDE):
 \be Y_t=\xi + \int_t^T g(s,Y_s,Z_s)ds-\int_t^T\langle Z_s,d
    B_s\rangle ,\q 0\leq t \leq T. \label{BSDE0}
 \ee
 Under conditions (C1)-(C3), it has a unique adapted solution
 $(Y,Z)$. In the following, denote $Y_t$ by
 $\cE_{t,T}^g[\xi]$ to emphasize the dependence on the
 parameter $(\xi, g)$ and the initial and terminal times pair $(t,T)$.
 Coquet, Hu, Memin and Peng~\cite{CHMP} introduce the notion of an $\{\cF_t, 0\le t\le T\}$-consistent expectation, which is defined to be a nonlinear
 functional $\cE[\cdot]$ on $L^2(\cF_T)$,
 satisfying the following three axioms (E1)-(E3).

{\rm  (E1) } $\cE[Y_1]\le \cE[Y_2]$ if $Y_1,Y_2\in L^2(\cF_T)$ and
$Y_1\le Y_2$. Furthermore, if $Y_1\le Y_2$, then $Y_1=Y_2$ if
 $\cE[Y_1]=\cE[Y_2]$.

{\rm  (E2) } $\cE[C]=C$ for any constant $C$. And

{\rm (E3) } The $\cE$-expectation (denoted by $\cE[Y|\cF_t]$)  of
$Y\in L^2(\cF_T)$ conditioned on $\cF_t$ exists uniquely for any
$t\in
 [0,T]$. That is, there is unique $\eta\in L^2(\cF_t)$ such that
 \begin{equation}
\cE[\eta 1_A]=\cE[Y 1_A],\q \forall A\in \cF_t.
\end{equation}
The third axiom (E3) shows that an $\{\cF_t, 0\le t\le
T\}$-consistent expectation $\cE[\cdot]$ induces an $\{\cF_t, 0\le
t\le T\}$-consistent dynamic operator, which will be denoted by
$\{\cE_{s,t}[\cdot], 0\le s\le t\le T\}$. We shall identify one of
them  with the other.

Coquet, Hu, Memin and Peng~\cite{CHMP} proved the
 following representation result.

 \bt \label{BSDEConverse} \it Assume that an $\{\cF_t, 0\le t\le T\}$-consistent
 expectation $\{\cE_{s,t}, 0\le s\le t\le T\}$  satisfies (H1)' for some positive constant $\mu$ and the following assumption:

\rm (H2)' \it  {\bf Translation invariance.}
$\mathcal{E}_{s,t}[X+Y]=\mathcal{E}_{s,t}[X]+Y\ \as $ for $s\le t,
X\in L^2(\cF_t)$ and $Y\in L^2(\cF_s)$.

\noindent Then, there is a random field $g: \O\times [0, T]\times
\mathbb{R}^d\to \mathbb{R}$
 such that

 (i) $g(t,0)=0$ for $\ae t\in [0,T]$ and $g(\cdot, z)\in \cL^2_\cF(0,T)$ for $z\in \mathbb{R}^d$,

 (ii) $|g(t,z_1)-g(t,z_2)|\le \mu |z_1-z_2|$ for $z_1, z_2\in \mathbb{R}^d$ and $\ae t\in [0,T]$, and

 (iii) $\cE_{0,T}[Y]=\cE^g_{0,T}[Y]$ for any $Y\in L^2(\cF_T)$. \et\rm

 Note that an $\{\cF_t, 0\le t\le T\}$-consistent expectation $\{\cE_{s,t}, 0\le
s\le t\le T\}$  has the following four properties.

 (A1) {\bf Strict Monotonicity.} $\mathcal{E}_{s,t}[Y]\geq
\mathcal{E}_{s,t}[Y^{'}]$ if $Y,Y'\in L^2(\cF_t)$ and $Y\geq
Y^{'}$. Furthermore, if $Y'\ge Y$ and
$\cE_{s,t}[Y']=\cE_{s,t}[Y]$, then $Y'=Y$.

(A2) $\mathcal{E}_{s,t}[Y]=Y \, \as $ for $s\in [0,t]$ and $Y\in
L^2(\cF_s)$.

(A3) {\bf Time consistency.}
$\mathcal{E}_{r,s}[\mathcal{E}_{s,t}[Y]]=\mathcal{E}_{r,t}[Y]$ if
$r\leq s\leq t\leq T$ and $Y\in L^2(\cF_t)$.

(A4) {\bf Zero-One law.}
$\mathcal{E}_{s,t}[1_{A}Y]=1_{A}\mathcal{E}_{s,t}[Y]$\  a.s. for
$s\leq t, A\in \mathcal{F}_s, $ and $Y\in L^2(\cF_t)$.\\
Theorem~\ref{BSDEConverse} means that an $\{\cF_t, 0\le t\le
T\}$-consistent expectation always has a BSDE representation if it
is both $\cE^\mu$-dominated for some $\mu>0$ and
constant-preserving. Note that (H1) is stronger than (H1)' and
however, the latter together with (H2)' and (A1)-(A4) implies the
former (see Remark~\ref{h1h1'} for further details). Further,
(H2)' is stronger than (H2). Therefore, Theorem~\ref{BSDEConverse}
is the special case of $S\equiv -\infty$ for
Theorem~\ref{RBSDEConverse}.

Though Theorem~\ref{RBSDEConverse} is a generalized version of
Theorem~\ref{BSDEConverse}, it turns out to be far from
straightforward in both its formulation and its proof, due to the
appearance of the floor. Indeed, they call for original ideas. In
particular, two key points are worth to be mentioned here. The one
is to make full use of the non-increasing and
floor-above-invariant assumption (H2) of forward translation and
the assumption (H1) of $\cE^\mu$-super-domination for some
$\mu>0$, to extend the underlying $\{\cF_t, 0\le t\le
T\}$-consistent dynamic operator from the subset $L^2(\cF_T;S_T)$
of floor-dominating square-integrable random variables to the
whole space $L^2(\cF_T)$ of square-integrable random variables.
The extended dynamic operators are shown to be identified to an
$\{\cF_t, 0\le t\le T\}$-consistent expectation $\widetilde\cE$.
The generator $g$ of its BSDE representation given by
Theorem~\ref{BSDEConverse} turns out to be that of the desired
RBSDE. More interesting is the second key point---a beautiful
observation: the process $\{\cE_{t,T}[X], 0\le t\le T\}$ turns out
to be an $\widetilde \cE$-supermartingale for each $X\in
L^2(\cF_T;S_T)$. This fact allows us to apply a nonlinear
Doob-Meyer's decomposition theorem in~\cite{CHMP}, to give the
increasing process of $\{\cE_{t,T}[X], 0\le t\le T\}$ as a
solution to some BSDE reflected upwards on the floor. Eventually,
our proof of Theorem~\ref{RBSDEConverse} is both natural and
elegant. To diverse the difficulty in the above arguments, the
whole proof is divided into three sections: Sections 2-4. The
above two key points are exposed in detail separately in the first
two sections: Sections 2 and 3.

The rest of the paper is organized as follows. In Section 2, we
consider the case of the zero floor. In this case, for any
nonnegative terminal value, the process $\{K_t, 0\le t\le T\}$ in
RBSDE~(\ref{RBSDE0}) is zero and thus the RBSDE is reduced to a
BSDE. Nonlinear Doob-Meyer's decomposition theorems are not needed
in this case. We concentrate our attention to show how to extend
the family of dynamic operators defined on the subset
$L^2_+(\cF_T)$ of $L^2(\cF_T)$ to an $\{\cF_t, 0\le t\le
T\}$-expectation, which is defined on the whole space
$L^2(\cF_T)$.  In Section 3, we consider the case of the negative
floor. Restricting the underlying $\{\cE_{t,T}, 0\le t\le
T\}$-consistent dynamic operator to $L^2_+(\cF_T)$, we get an
$\{\cF_t, 0\le t\le T\}$-expectation by extending the restriction
to $L^2(\cF_T)$  in the way as shown in the preceding section. The
$\{\cF_t, 0\le t\le T\}$-expectation gives the generator $g$ of an
RBSDE by Theorem~\ref{BSDEConverse}. In addition, we have to give
the amount $\{K_t^X, 0\le t\le T\}$ to push upwards for
$\{\cE_{t,T}[X], 0\le t\le T\}$ with $X\in L^2(\cF_T;S_T)$, that
is, the increasing process in relevant to an RBSDE. For this
purpose, a nonlinear Doob-Meyer's decomposition theorems is shown
how to be used. The key point is to observe that the process
$\{\cE_{t,T}[X], 0\le t\le T\}$ turns out to be an $\widetilde
\cE$-supermartingale for each $X\in L^2(\cF_T;S_T)$. In Section 4,
we consider the general case of the upper bounded floor and give
the proof of Theorem~\ref{RBSDEConverse}.

\br The upper-bounded assumption on the floor $S$ in (H3) can be
removed. This subject will be detailed elsewhere. \er

\section{The case of the zero floor. }

In the case of the zero floor, the  properties (D2) and (D4)  read

(D2)' $\mathcal{E}_{s,t}[Y]=Y$ a.s. for $Y\in
L_{+}^2(\mathcal{F}_s)$ and $s\in [0,t]$.

(D4)' { \bf Zero-one law }
$\mathcal{E}_{s,t}[1_{A}Y]=1_{A}\mathcal{E}_{s,t}[Y]$, a.s. for
$s\in [0,t]$, $A\in \mathcal{F}_s$,  and $Y\in
L_{+}^2(\mathcal{F}_t)$.

From (D1), (D2)', (H1) and  (H2), we can give the following:

(H2)' We have for $ t\in [0,T], X \in L_{+}^2(\mathcal{F}_T)$ and
$Y\in L_{+}^2(\mathcal{F}_t)$,
$$\
\mathcal{E}_{t,T}[X+Y]=\mathcal{E}_{t,T}[X]+Y \, \as.$$

In fact, consider $ X \in L_{+}^2(\mathcal{F}_T)$ and $Y\in
L_{+}^2(\mathcal{F}_t)$. It follows from (D2)' that
$\cE_{t,T}[\epsilon]=\epsilon$ for $t\in [0,T]$ and very constant
$\epsilon>0$. Therefore, from (D1), we deduce that
$\cE_{t,T}[\epsilon+X]\ge \cE_{t,T}[\epsilon]=\epsilon>0\ \as $
for $t\in [0,T]$. Consequently, by (H2), we conclude that
\be\cE_{t,T}[\epsilon+X+Y]= \cE_{t,T}[\epsilon+X]+Y\ \as. \ee By
(H1), we can take the limit $\epsilon\to 0$, and we have the
desired equality.

It is easy to prove the following theorem:

\bt ~\label{kt0} Let $(y,z)$ be the adapted solution of
BSDE~(\ref{BSDE0}) with the terminal condition $Y\in\lft$. Then
the triple $(y,z,0)$ is the adapted solution of
RBSDE~(\ref{RBSDE0}) with the parameter $(Y,g,0)$. \et

\br Theorem~\ref{kt0} shows that the adapted solution of
RBSDE~(\ref{RBSDE0}) is the one of the corresponding
BSDE~(\ref{BSDE0}) for any terminal value $Y\in
L_{+}^2(\mathcal{F}_T)$ since the obstacle $\{S_t, 0\le t\le T\}$
is not active in this case.
 \er

 The following lemma is immediate and will be used later.

\bl \label{bsdeEst} We have

(i)  $E\bigl[\cE^\mu_{t,T}[X]^p\bigr]\le
\exp(2^{-1}p(p-1)^{-1}\mu^2(T-t)) E[X^p]$ for all $\mu>0, t\in
[0,T]$, and $X\in L^p(\cF_T)$, with $p\in (1,2];$

(ii) $\cE^{-\mu}_{t,T}[X+Y]=Y-\cE^\mu_{t,T}[-X] \, \as $ for  all
$\mu>0, t\in [0,T], X\in L^2(\cF_T)$, and $Y\in L^2(\cF_t);$ And

(iii)  $\cE^{\mu}_{0,T}[\cdot]$ and $\cE^{-\mu}_{0,T}[\cdot]$ are
strongly continuous in $L^2(\cF_T)$.\el

{\bf Proof. } We only prove the first assertion (i). The other two
assertions are easy to see.

For simplicity of notations, set $y_t:=\cE^\mu_{t,T}[X].$ Then, be
definition, there is unique $z\in \cL^2_\cF(0,T;\mathbb{R}^d)$
such that $(y,z)$ is the unique adapted solution of
BSDE~(\ref{bsde0}). Using It\^o's formula, we have \be
E|y_t|^p+{1\over 2} p(p-1)\int_t^T|y_s|^{p-2}|z_s|^2\, ds
=E|X|^p+2\mu p\int_t^T|y_s|^{p/2}|y_s|^{p/2-1}|z_s|\, ds. \ee
Since \be \mu p|y_s|^{p/2}|y_s|^{p/2-1}|z_s|\le {p\over
2(p-1)}\mu^2|y_s|^p+{1\over 2} p(p-1)|y_s|^{p-2}|z_s|^2,\ee we
have \be E|y_t|^p\le E|X|^p+{p\over 2(p-1)}\mu^2 \int_t^T|y_s|^p\,
ds. \ee The standard arguments of using Gronwall's inequality then
gives the desired inequality. \endpf

\br See Coquet, Hu, Memin and
 Peng~\cite[Lemma 2.2, page 5]{CHMP} for the detailed proof of the first assertion in the case of $p=2$, which is easy and
 standard.\er

In view of  Theorem~\ref{kt0}, we have from Peng's work~\cite{P2}
on the properties of BSDE the following

\bt \label{thmnew} Assume that conditions (C2)-(C4) are satisfied.
Moreover, assume that $g(\cdot,\cdot,0)\equiv 0$. Then,
$\{\cE^{r,g}_{s,t}, 0\le s\le t\le T\}$ satisfies (H1), (D1),
(D2)', (D3), (D4)' and (H2)'. Therefore, it is an $\{\cF_t, 0\le
t\le T\}$-consistent dynamic operator with the zero floor. \et

Roughly speaking, Theorem~\ref{thmnew} asserts that an RBSDE with
the obstacle being zero introduces an $\{\cF_t, 0\le t\le
T\}$-consistent dynamic operator.  In what follows, we shall
consider the converse problem. That is, we shall associate an
$\{\cF_t, 0\le t\le T\}$-consistent dynamic operator with the zero
floor to a BSDE reflected on the zero
 floor. For this purpose, we  establish the following six
 preliminary lemmas. First, we introduce some notations.

 \bde  For an $\{\cF_t,
0\le t\le T\}$-consistent dynamic operator $\{\cE_{s,t}[\cdot],
0\le s\le t\le T\}$, define the system of dynamic operators
$\mathcal{E}[\cdot |\mathcal{F}_t]: L_{+}^2(\mathcal{F}_T)
 \longrightarrow L_{+}^2(\mathcal{F}_t), 0\leq t<
 T$ by
\begin{equation}
\mathcal{E}[Y|\mathcal{F}_t]:=\mathcal{E}_{t,T}[Y], a.s. \hbox{
\rm for } Y\in L_{+}^2(\mathcal{F}_T)
\end{equation}
and the nonlinear functional $\cE[\cdot]: L^2(\cF_T)\to
\mathbb{R}$ by
\begin{equation}
\mathcal{E}[Y]:=\mathcal{E}_{0,T}[Y] \hbox{ \rm for }  Y\in
L_{+}^2(\mathcal{F}_T).
 \end{equation}
 \ede

The two notations $\cE[\cdot]$ and $\mathcal{E}[\cdot
|\mathcal{F}_t]$ behave in $L^2_+(\cF_T)$ exactly like an
$\{\cF_t, 0\le t\le T\}$-consistent expectation and its
conditional $\{\cF_t, 0\le t\le T\}$-consistent expectation on
$\cF_t$. The only difference lies in the domains of variables: the
former's are $L^2_+(\cF_T)$, while the latter's are $L^2(\cF_T)$.
This can be seen from the following obvious lemma.

\bl  For an $\{\cF_t, 0\le t\le T\}$-consistent  dynamic operator
$\{\cE_{s,t}[\cdot], 0\le s\le t\le T\}$, $\cE[\cdot]$ and
$\mathcal{E}[\cdot |\mathcal{F}_t]$ have the following properties
(E1)'-(E3)'.

{\rm  (E1)' } $\cE[Y_1]\le \cE[Y_2]$ if $Y_1,Y_2\in L^2_+(\cF_T)$
and $Y_1\le Y_2$. Furthermore, if $Y_1\le Y_2$, then $Y_1=Y_2$ if
 $\cE[Y_1]=\cE[Y_2]$.

{\rm  (E2)' } $\cE[C]=C$ for any constant $\widetilde C>0$. And

{\rm (E3)' } For any $Y\in L^2_+(\cF_T)$, there is unique $\eta\in
L^2_+(\cF_t)$ such that
 \begin{equation}
\cE[\eta 1_A]=\cE[Y 1_A],\q \forall A\in \cF_t,
\end{equation} which is equal to $\cE[Y|\cF_t]$. Therefore, $\cE[Y|\cF_t]$
can be viewed as the $\cE$-expectation  of $Y\in L^2_+(\cF_T)$
conditioned on $\cF_t$ for any $t\in
 [0,T]$, though $\cE$ is in general not an $\{\cF_t, 0\le t\le T\}$-consistent expectation at all.  \el

  \bl~\label{conti} Assume that  $\e[Y|\f_t]$ satisfies
 conditions (H1)' for $\mu> 0$. Then,
 \begin{equation}
 \ |\e[\xi_1]-\e[\xi_2]|\leq \exp({{1\over
2}{\mu}^2T})||
 \xi_1-\xi_2||_{L^2}, \q \forall \xi_1, \xi_2 \in L_{+}^2(\f_T).
 \end{equation} Therefore,
$\e[\cdot]$ is strongly continuous in $L^2(\f_T)$.
 \el

 {\bf Proof. }   Using (D1) and (H1)', we have
 \begin{eqnarray}
 \e[\xi_1]-\e[\xi_2] &\leq& \e[|\xi_1-\xi_2|+\xi_2]-\e[\xi_2] \nonumber \\
  &\leq& \e^{\mu}[|\xi_1-\xi_2|]. \nonumber
 \end{eqnarray}
 Using Lemma~\ref{bsdeEst}, we have
 \begin{equation}~\label{inequmu}
 \e^{\mu}[|\xi_1-\xi_2|]^2\leq \exp{({\mu}^2T)}E[|\xi_1-\xi_2|^2].
 \end{equation}
 Therefore,
 \begin{equation}
\e[\xi_1]-\e[\xi_2]\leq \exp({{1\over
2}{\mu}^2T})(E|\xi_1-\xi_2|^2)^{1/2}=\exp({{1\over
2}{\mu}^2T})||\xi_1-\xi_2||_{L^2},\q \forall \xi_1, \xi_2 \in
L_{+}^2(\f_T).
 \end{equation}

 Identically, we can show
  \begin{equation}
\e[\xi_2]-\e[\xi_1]\leq \exp({{1\over
2}{\mu}^2T})||\xi_1-\xi_2||_{L^2},\q \forall \xi_1, \xi_2 \in
L_{+}^2(\f_T).
 \end{equation}
 The proof is then complete. \endpf

\bl\label{binequ} ({\bf $\e^{\mu}$-domination}) Assume that
$\e[Y|\f_t]$ satisfies (H1)' and (H2)'. Then, we have
\begin{equation}
 \e^{- \mu}[Y]\leq \e[X+Y]-\e[X] \leq \e^{\mu}[Y], \q \forall
X,Y\in L_{+}^2(\f_T).
\end{equation}
\el

{\bf Proof. } In view of (H1)', it is sufficient to prove the
following \be \label{des} \cE[X+Y]-\cE[Y]\ge \cE^{-\mu}[Y], \q
\forall X,Y\in L^2_+(\cF_T).\ee

If $Y\leq n \ \as$ for some integer $n$, then using (H2)', (H1)',
and Lemma~\ref{bsdeEst}, we have
\begin{eqnarray}
n-(\e[X+Y]-\e[X]) &=& n+\cE[X]-\cE[X+Y]=\e[X+n]-\e[X+Y]\nonumber \\
                &\leq & \e^{\mu}[n-Y]\le n-\cE^{-\mu}[Y]. \nonumber
\end{eqnarray}
Therefore,
\begin{equation}~\label{tem} \e[X+Y]-\e[X]\ge \e^{- \mu}[Y].
\end{equation}

 In general,  consider $Y\in L_{+}^2(\f_T)$.  Define
$Y_n:=Y1_{\{Y\leq n\}}.$ Then, $Y_n$ converges to $Y$ strongly in
$L^2(\cF_T)$.  The above arguments show that
 $$ \e[X+Y_n]-\e[X] \ge \e^{-
\mu}[Y_n].
$$
Then using Lemma~\ref{conti}, the desired result~(\ref{des})
follows by passing to the limit $n\to \infty$ in the last
inequality.\endpf

\br The above proof of Lemma~\ref{binequ} is more complicated
than~\cite{CHMP} since both assumptions (H1)' and (H2)' are
weaker.\er

For $\zeta\in L_{+}^2(\f_T)$, define the operator
$\e^{\zeta}_{s,t}[\cdot]: L_{+}^2(\f_t)\to L^2(\cF_s)$ by
\begin{equation}
\e^{\zeta}_{s,t}[X]:=\e_{s,t}[X+\zeta]-\e_{s,t}[\zeta], \q \forall
X\in L^2_+(\cF_T).
\end{equation}

Identically as in the case of an $\{\cF_t, 0\le t\le
T\}$-expectation (see~\cite[Lemma 4.3, page 10]{CHMP} for detailed
proof), we can show the following lemma.

\bl Let $\zeta\in L^2_+(\cF_T)$. If $\{\e_{s,t}[\cdot], 0\le s\le
t\le T\}$ is an $\{\cF_t, 0\le t\le T\}$-consistent dynamic
operator defined on $ L_{+}^2(\f_T)$ and satisfies (H1)' and
(H2)', then the operator $\{\e^{\zeta}_{s,t}[\cdot], 0\le s\le
t\le T\}$ is also an $\{\cF_t, 0\le t\le T\}$-consistent dynamic
operator defined on $ L_{+}^2(\f_T)$, and satisfy (H1)' and (H2)'.
The expectation $\e^{\zeta}[X|\f_t]$ of $X\in L^2_+(\cF_T)$
conditioned on $\f_t$ is given by the formula:
\begin{equation}
\e^{\zeta}[X|\f_t]=\e[X+\zeta |\f_t]-\e[\zeta |\f_t].
\end{equation}
\el

\bl~\label{ineq} Assume that the two $\f$-consistent dynamic
operators $\e^1[\cdot]$ and $\e^2[\cdot]$ defined on $
L_{+}^2(\f_T)$ satisfy (H1)' and (H2)'. If
$$\e^1[X]\leq \e^2[X], \q \forall X\in \lft,$$
then for all t,
$$\e^1[X|\f_t]\leq \e^2[X|\f_t],\q \as \hbox{ \rm for all } X\in \lft. $$
\el

{\bf Proof.} The proof is divided into the two steps.

\bf Step 1. \rm The case of  $ X\leq n$ for some positive integer
$n$. Set
$$\eta=\e^2[X|\f_t]-\e^1[X|\f_t].$$ Then
$$-\eta 1_{\{\eta \leq 0\}}=\e^1[X1_{\{\eta \leq 0\}}|\f_t]-
\e^2[X1_{\{\eta \leq 0\}}|\f_t]\geq 0
$$ and
\begin{eqnarray*}
n &\leq& \e^1[-\eta1_{\{\eta \leq 0\}}+n] \nonumber \\
  &=& \e^1[\e^1[X1_{\{\eta \leq 0\}}|\f_t]-
\e^2[X1_{\{\eta \leq 0\}}|\f_t]+n] \nonumber \\
&=& \e^1[X1_{\{\eta \leq 0\}}-
\e^2[X1_{\{\eta \leq 0\}}|\f_t]+n] \nonumber \\
&\leq& \e^2[X1_{\{\eta \leq 0\}}-
\e^2[X1_{\{\eta \leq 0\}}|\f_t]+n] \nonumber \\
&=& \e^2[\e^2[X1_{\{\eta \leq 0\}}-
\e^2[X1_{\{\eta \leq 0\}}|\f_t]+n|\f_t]] \nonumber \\
&=& \e^2[\e^2[X1_{\{\eta \leq 0\}}|\f_t]+n-
\e^2[X1_{\{\eta \leq 0\}}|\f_t]] \nonumber \\
&=& n.
\end{eqnarray*}
Thus,
$$\e^1[-\eta1_{\{\eta \leq 0\}}+n]=n.$$
From the floor-above  strict monotonicity (D1), we have
$$-\eta1_{\{\eta \leq
0\}}=0.
$$ That is,  $P(\{\eta \leq 0\})=0$, \q a.s.. Thus,
$$\e^1[X|\f_t]\leq \e^2[X|\f_t],\q a.s.. $$

\bf Step 2. \rm The general case of  $ X\in \lft$.

For $ X\in \lft$,  define the truncation $X_n:=X1_{\{X\leq n\}}$.
Then it follows from Lemma~\ref{conti} that
$$\lim_{n\rightarrow \infty}\e^1[X_n]=\e^1[X].$$
From Step 1, we have
$$\e^1[X_n|\f_t]\leq \e^2[X_n|\f_t], \q n=1,2,\ldots.$$ Obviously, it is sufficient to prove the following
\begin{equation}~\label{3}
\lim_{n\rightarrow \infty}\e^1[X_n|\f_t]=\e^1[X|\f_t],\ \forall
X\in \lft
\end{equation}
We  now prove it by contradiction.

Otherwise, there exists $ 0\leq \epsilon \leq 1$ and $A\in \f_t$
such that $P(A)> 0$ and
$$\e^1[X_n|\f_t]1_A\leq (\e^1[X|\f_t]-\epsilon)1_A. $$
Since $\cE^1[X_n+1|\cF_t]=\cE^1[X_n|\cF_t]+1$ and
$\cE^1[X+1|\cF_t]=\cE^1[X|\cF_t]+1$ (by (H2)'), we have
$$\e^1[X_n+1|\f_t]1_A\leq (\e^1[X+1|\f_t]-\epsilon)1_A. $$
Using (D1), we have $\e^1[\e^1[X_n+1|\f_t]1_A] \leq
\e^1[(\e^1[X+1|\f_t]-\epsilon)1_A]$.  Then letting $n\to \infty,$
we have
\begin{eqnarray}~\label{1}
\lim_{n\rightarrow \infty}\e^1[\e^1[X_n+1|\f_t]1_A] &\leq&
\lim_{n\rightarrow
\infty}\e^1[(\e^1[X+1|\f_t]-\epsilon)1_A] \nonumber \\
&<& \e^1[(X+1)1_A].
\end{eqnarray}
While we have by Lemma~\ref{bsdeEst} the following
\begin{eqnarray}~\label{2}
\lim_{n\rightarrow \infty}\e^1[\e^1[X_n+1|\f_t]1_A] &=&
\lim_{n\rightarrow
\infty}\e^1[(X_n+1)1_A] \nonumber \\
&=& \e^1[(X+1)1_A].
\end{eqnarray}
This contradicts (\ref{1}). Therefore, (\ref{3}) is true. \endpf

Combining Lemmas~\ref{ineq} and~\ref{binequ}, we have

\bl~\label{inquft} Let $\{\e_{s,t}[\cdot], 0\le s\le t\le T\}$ be
an $\f$-consistent dynamic operator defined on $\lft$ and satisfy
(H1)' and (H2)'. Then, for each $t\leq T$, we have
\begin{equation}
\e^{- \mu}[Y|\f_t]\leq \e[X+Y|\f_t]-\e[X|\f_t] \leq
\e^{\mu}[Y|\f_t], \as \hbox{ \rm for all }  X,Y\in L_{+}^2(\f_T).
\end{equation}
\el

\br \label{h1h1'} Lemma~\ref{inquft} shows that (H1)' together
with (H2)' and (D1)-(D4) implies (H1), as pointed in the
introduction.\er

 In the following, we extend
the underlying $\{\cF_t, 0\le t\le T\}$-consistent dynamic
operator from the subset $L^2_+(\cF_T)$ of nonnegative
square-integrable random variables to the whole space $L^2(\cF_T)$
of square-integrable random variables. As a preliminary, we have
the following fact.

 \bl \label{lmm2.7} Let  $\{\e_{s,t}[\cdot], 0\le s\le t\le T\}$ be
an $\f$-consistent dynamic operator defined on $\lft$ and satisfy
(H1)' and (H2)'. Let $X_n:=X1_{\{X\geq -n\}}$ and
$Y_n:=\e[X_n+n|\f_t]-n$ for $\, X\in \lf$ and $n=1,2,\ldots$. Then
$\{Y_n\}_{n=1}^\infty$
 is a Cauchy sequence in $L^2(\f_t)$ equipped with the norm
 $||\cdot||=\sqrt{E|\cdot|^2}$. If $X\in L^2_+(\cF_T)$, then
 $X_n=X$ and $Y_n=\cE[X|\cF_t]$ for $n=1,2,\ldots$.
\el

{\bf Proof.}  For the two positive integers $m$ and $n$ such that
$m>n$, we have
\begin{eqnarray}
Y_m-Y_n &=& \e[X_m+m|\f_t]-m-(\e[X_n+n|\f_t]-n) \nonumber \\
&=& \e[X_m+m|\f_t]-\e[X_n+m|\f_t] \nonumber \\
&=& \e[X_n+X1_{\{-m\leq X\leq -n\}}+m|\f_t]-\e[X_n+m|\f_t].
\end{eqnarray}
Thus, from Lemmas~\ref{inquft} and \ref{bsdeEst}, we have
\begin{eqnarray}
E(Y_m-Y_n)^2 &=& E(\e[X_n+m|\f_t]-\e[X_n+X1_{\{-m\leq X\leq -n\}}+m|\f_t])^2 \nonumber \\
        &\leq & E(\e^{\mu}[-X1_{\{-m\leq X\leq -n\}}|\f_t])^2 \nonumber \\
        &\leq & e^{{\mu}^2(T-t)}E[X^21_{\{-m\leq X\leq -n\}}].
\end{eqnarray}
Since $X\in L^2(\cF_T)$, we have
$$\lim_{m,n\to \infty}E[X^21_{\{-m\leq X\leq -n\}}]=0.
$$
Therefore,
$$\lim_{n,m\to \infty}E(Y_m-Y_n)^2=0.$$
\endpf

Lemma~\ref{lmm2.7} shows that $\widehat \e[X|\f_t]$ introduced
below is well defined for any $\ X\in \lf$.

\bde~\label{expand} For  $X\in \lf$, denote
\begin{equation}
X_n:=X1_{\{X\geq -n\}}.
\end{equation}
We introduce a dynamic operator $\{\widehat \cE_{s,t}[X], X\in
L^2(\cF_T); 0\le s\le t\le T\}$ by
\begin{equation}
\widehat \e[X_n|\f_t]:=\e[X_n+n|\f_t]-n, \q n=1,2,\ldots;
\end{equation}
and
\begin{equation}~\label{expansion}
 \widehat\e[X|\f_t]:=\lim_{n\rightarrow
\infty}\widehat\e[X_n|\f_t]
\end{equation}
for any  for $ X\in \lf$. \ede

\br \label{iden}It is obvious that if $X\in L^2_+(\cF_T)$, we have
$\widehat \cE[X|\cF_t]=\cE[X|\cF_t]\, \as$ for $t\in [0,T]$.\er

We have the following continuity on the extended operator
$\widehat \cE[\cdot|\cF_t]$.

\bt ~\label{lim} Let $\{\e_{s,t}[\cdot], 0\le s\le t\le T\}$ be an
$\f$-consistent dynamic operator defined on $\lft$ and satisfy
(H1)' and (H2)'.  For each $t\in [0,T]$, the conditional
expectation operator $\widehat \cE[\cdot|\cF_t]$ is strongly
continuous from $L^2_+(\cF_T)$ to $L^2(\cF_t)$. That is, if
$$\lim_{n\to\infty}X_n=X, \hbox{ \rm strongly in } \lf,
$$then
\begin{equation}
\lim_{n\to \infty}\widehat\e[X_n|\f_t]=\widehat \e[X|\f_t],\q
\hbox{ \rm strongly in } L^2(\cF_t).
\end{equation}
\et

{\bf Proof.} For the two positive integers $m$ and $n$ such that
$m> n$, define
$$\d_1:=\widehat\e[X_n|\f_t]-\e[X_n1_{\{X_n\geq -m\}}+m|\f_t]-m$$
and
$$\d_2:=(\e[X_n1_{\{X_n\geq
-m\}}+m|\f_t]-m)-(\e[X1_{\{X\geq -m\}}+m|\f_t]-m).$$
 Using Lemma~\ref{inquft}, we have
\begin{eqnarray*}
|\triangle_2| &\leq &  |\e[X 1_{\{X\geq -m\}} +m +|X_n
1_{\{X_n\geq -m\}}-X 1_{\{X\geq -m\}}||\f_t] \nonumber \\ & &
-\e[X1_{\{X\geq -m\}}+m|\f_t]|
\nonumber \\
&\leq &  \e^{\mu}[|X_n 1_{\{X_n\geq -m\}}-X 1_{\{X\geq
-m\}}||\f_t]\nonumber \\
&\le& \e^{\mu}[|X_n 1_{\{X_n\geq
-m\}}-X_n||\f_t]+\e^{\mu}[|X_n-X||\f_t]+\e^{\mu}[|X-X 1_{\{X\geq -m\}}||\f_t]. \\
\end{eqnarray*}
Thus, using Lemma~\ref{bsdeEst}, we have
\begin{eqnarray*}
E|\d_2|^2 &\leq& 3E\e^{\mu}[|X_n 1_{\{X_n\geq
-m\}}-X_n||\f_t]^2+3E\e^{\mu}[|X_n-X||\f_t]^2\nonumber \\
& &+3E\e^{\mu}[|X-X 1_{\{X\geq -m\}}||\f_t]^2\nonumber \\
&\leq& 3e^{{\mu}^2(T-t)}E|X_n 1_{\{X_n\geq
-m\}}-X_n|^2\nonumber \\
& &+ 3e^{{\mu}^2(T-t)}(E|X_n-X|^2
+E|X-X 1_{\{X\geq -m\}}|^2 )\nonumber \\
&\leq& 3e^{{\mu}^2T}E|X_n 1_{\{X_n< -m\}}|^2
+3e^{{\mu}^2T}(E|X_n-X|^2+E|X 1_{\{X< -m\}}|^2 ). \nonumber
\end{eqnarray*}
Letting $m\to \infty$ in the following
\begin{eqnarray*}
&&E|\widehat\e[X_n|\f_t]-\widehat\e[X|\f_t]|^2 \\
&\leq&
3(E|\d_1|^2+E|\d_2|^2+E|\e[X1_{\{X\geq
-m\}}+m|\f_t]  -m-\widehat\e[X|\f_t]|^2)\nonumber \\
&\leq & C_1\biggl(E\d_1^2+E|X_n 1_{\{X_n< -m\}}|^2+E|X_n-X|^2+E|X
1_{\{X< -m\}}|^2\nonumber
\\&& +E|\e[X1_{\{X\geq -m\}}+m|\f_t]
-m-\widehat\e[X|\f_t]|^2\biggr),
\end{eqnarray*}
we have
\begin{eqnarray*}
E|\widehat\e[X_n|\f_t]-\widehat\e[X|\f_t]|^2 \leq  C_1E|X_n-X|^2.
\end{eqnarray*}
Therefore,
$$\lim_{n\rightarrow \infty}\widehat\e[X_n|\f_t]=\widehat\e[X|\f_t].
$$
\endpf

The extended dynamic operator $\widehat \cE$  can be shown to be
identified to an $\{\cF_t, 0\le t\le T\}$-consistent expectation
$\widehat\cE$.

 \bt \label{thm2.4} Let $\{\e_{s,t}[\cdot], 0\le s\le t\le T\}$ be an
$\{\f_t, 0\le t\le T\}$-consistent dynamic operator defined on
$\lft$ and satisfy (H1)' and (H2)'. Then, the  dynamic operator
$\widehat\e[\cdot|\f_t]$ defined by (\ref{expansion}) satisfies
(H1)', (H2)', (D1), (D2)', (D3), and (D4)'. Therefore, it is an
$\{\cF_t, 0\le t\le T\}$-consistent expectation. \et

{\bf Proof.} For $X\in \lf$ and $Y\in \lf$, define
$$X_n:=X1_{\{X\geq -n\}},\ Y_n=Y1_{\{Y\geq -n\}}.$$
It is easy to check all the properties of Monotonicity (D1),
Lemma~\ref{inequmu}, (D2)', (D3), (D4)', (H2)' and (H1)'  for
$X_n$ and $Y_n$. From Theorem~\ref{lim}, let $n\to \infty$,
monotonicity (D1), Lemma~\ref{inequmu}, (D2)', (D3), (D4)', (H2)',
and (H1)' then follow immediately.

If $Y\geq X$ and $\widehat\e[Y|\f_t]=\widehat\e[X|\f_t]$, then it
follows from Lemma~\ref{inequmu} that
 $$0=\widehat\e[Y|\f_t]-\widehat\e[X|\f_t]\geq \e^{-\mu}[Y-X|\f_t]
\geq 0.$$ Thus,$$\e^{-\mu}[Y-X|\f_t] =0.$$ Therefore,$$Y=X.$$ This
shows the strict monotonicity.
\endpf

The generator $g$ of the BSDE representation of the $\{\cF_t, 0\le
t\le T\}$-consistent expectation $\widehat\e[\cdot]$ given by
Theorem~\ref{BSDEConverse} turns out to be that of the desired
RBSDE.

\bt ~\label{existg} Theorem~\ref{RBSDEConverse} is true in the
case of the zero floor  $S\equiv 0$. Moreover, the assumption (H1)
may be weakened to (H1)'. \et

{\bf Proof.} Theorem~\ref{thm2.4} shows that $\widehat\e[X], X\in
\lf$,  defined by (\ref{expansion}) is an
 $\{\cF_t, 0\le t\le T\}$-expectation, and satisfy (H1)' and (H2)'. From Theorem~\ref{BSDEConverse}, there exists a
 function $g:\O\times [0,T]\times \mathbb{R}^d\to \mathbb{R}$, satisfying (C2) and (C3) and the following
 properties: $g(0,z)=0, |g(t,z)|\leq \mu |z|$ for $\ae t\in
 [0,T]$, and $\widehat \cE[X]=\cE^g[X]$ for any $X\in L^2(\cF_T)$.

In version of  Remark~\ref{iden} and Lemma~\ref{bsdeEst}, we have
 $$\cE[X]=\widehat\cE[X]=\cE^g[X]=\cE^{r,g}[X] \hbox{ \rm for any}  X\in \lft.$$
 Hence, for any $A\in \cF_t$ and $X\in \lft$, we have
 \be \cE[X1_A]=\cE^{r,g}[X1_A].
 \ee
 This shows that $\cE[X|\cF_t]=\cE^{r,g}[X|\cF_t]$ for any $X\in
 \lft$. That is, $\cE_{t,T}[X]=\cE^{r,g}_{t,T}[X]$ for any $X\in
 \lft$.  \endpf

 \section{The case of the negative floor. }

Consider the dynamic nonlinear operator with a negative floor
$S\le 0$.

In what follows, we discuss the property of the solution Y to a
RBSDE~(\ref{RBSDE0}) from a view point of operators.  Given $t\leq
T$ and $Y\in
 L^2(\mathcal{F}_t;S_t)$, consider the following RBSDE on the
 interval $[0,t]$:
 \begin{equation}~\label{RBSDE}
 Y_s=\xi + \int_s^t g(r,Y_r,Z_r)dr+ K_t-K_r-\int_s^t\langle Z_r,d
    B_r\rangle,\q  0\leq s \leq t.
    \end{equation}

Suppose that the generator $g$ satisfies conditions (C2)-(C4).

\bde
 Define,  for  $0\leq s\leq t<\infty$ and $Y\in
 L^2(\f_t;S_t)$,
 \begin{equation}
 \mathcal{E}_{s,t}^{r;g,S}[Y]:=y_s ,
 \end{equation}  where $(y, z,K)$ is the solution of
 RBSDE~(\ref{RBSDE0})
with the parameter $(Y,g,S).$ \ede

\bt~\label{the prop of RBSDE} Let the random field $g$ satisfy
(C2)-(C4). Then if $S_t\le 0\ \as $ for $t\in [0,T]$, the dynamic
operator $\{\mathcal{E}_{s,t}^{r;g,S}, 0\le s\le t\leq T\}$
satisfies Axioms (H1), (H2), and (D1)-(D4). \et

{\bf Proof.} Let us prove (H1) first.
   Let $X_s^1$ and $X_s^2$ be the solutions of RBSDE~(\ref{RBSDE}) with
the terminal conditions $X+Y$ and $X$, respectively.
$Y_s=\e^{\mu}_{s,T}[Y]$. Applying It$\hat{o}$'s formula to
$|(X_s^1-X_s^2-Y_s)^{+}|^2$, and taking the expectation, we have:
\begin{eqnarray*}
& & E|(X_s^1-X_s^2-Y_s)^+|^2+E\int_s^t1_{\{X_r^1-X_r^2\geq Y_r\}}
|Z_r^1-Z_r^2-Z_r|^2dr \nonumber \\
&\leq &
2E\int_s^t(X_r^1-X_r^2-Y_r)^+[f(r,Z_r^1)-f(r,Z_r^2)-\mu|Z_r|]dr
\nonumber \\
&+& 2E\int_s^t(X_r^1-X_r^2-Y_r)^+(dK_r^1-dK_r^2). \nonumber
\end{eqnarray*}
Since  $X_r^1-X_r^2-Y_r\leq X_r^1-S_r$, we have
\begin{eqnarray*}
\int_s^t(X_r^1-X_r^2-Y_r)^+(dK_r^1-dK_r^2) & =&
-\int_s^t(X_r^1-X_r^2-Y_r)^+dK_r^2 \leq  0. \nonumber
\end{eqnarray*}
It follows from the conditions (C2) and (C3) that
\begin{eqnarray*}
& & E|(X_s^1-X_s^2-Y_s)^+|^2+E\int_s^t1_{\{X_r^1-X_r^2\geq Y_r\}}
|Z_r^1-Z_r^2-Z_r|^2dr \nonumber \\
&\leq & 2E\int_s^t(X_r^1-X_r^2-Y_r)^+\mu[Z_r^1-Z_r^2-Z_r|]dr
\nonumber \\
&\leq & E\int_s^t1_{\{X_r^1-X_r^2\geq Y_r\}} |Z_r^1-Z_r^2-Z_r|^2dr
+E\int_s^t\mu^2[Z_r^1-Z_r^2-Z_r|]^2ds.
\end{eqnarray*}
Hence
$$
E|(X_s^1-X_s^2-Y_s)^+|^2\leq E\int_s^t\mu^2[Z_r^1-Z_r^2-Z_r|]^2ds,
$$
and from Gronwall's lemma, $(X_s^1-X_s^2-Y_s)^+=0, \ \ 0\leq s
\leq t$. Therefore,
$$\mathcal{E}_{s,t}^{r;g,S}[X+Y]-\mathcal{E}_{s,t}^{r;g,S}[X]\leq
\mathcal{E}_{s,t}^{\mu}[Y], \q \forall X\in
L^{2}(\mathcal{F}_t;S_t) \hbox{ \rm and }Y\in
L_{+}^2(\mathcal{F}_t).
$$

Similarly, we can show the following inequality in (H2):
$$\mathcal{E}_{s,t}^{r;g,S}[X+Y]\leq \mathcal{E}_{s,t}^{r;g,S}[X]+
Y, \q  \forall X\in L^{2}(\mathcal{F}_t;S_t),Y\in
L_{+}^2(\mathcal{F}_s).$$

 If
$\e_{s,t}^{r;g,S}[X]>S_s \, \as$ for $s\in[r,t]$, in version of
the following equality \be
\int_0^t[\cE_{s,t}^{r;g,S}[X]-S_s]dK_s=0,\ee we have
$$K_s\equiv K_0=0, a.s. \ \forall  s\in[r,t]. \mbox{ \rm That is }  K\equiv 0, a.s..$$
Thus the solution of the RBSDE with the terminal condition $Y\geq
X$ is actually  that of  a BSDE. It follows from the strict
monotonicity of BSDEs that
 $$ Y=X, a.s. \Longleftrightarrow
Y\geq X, a.s. \mbox{ \rm and } \mathcal{E}_{r,t}[Y]=
\mathcal{E}_{r,t}[X], a.s..$$ The proof of (D1)-(D4) can be found
in El Karoui and Quenez~\cite{ELKaQuenez} and Xu~\cite{Xu}.
\endpf

Let $\{\e_{s,t}[\cdot], 0\le s\le t\le T\}$ be an $\f$-consistent
dynamic operator with floor $S\le 0$ and satisfy (H1) and (H2).
Under suitable conditions, we can show that it can be represented
as a BSDE reflected upwards on a negative obstacle. For this
purpose, consider its restriction on $L^2_+(\cF_T)$ and denote it
as
\begin{equation}
\e_{t,T}^{+}[X]:=\e_{t,T}[X], \q \forall X\in \lft.
\end{equation}
It is straightforward to check that $\e_{t,T}^{+}$ satisfies
 conditions (H1)', (H2)', (D1), (D2)', (D3) and (D4).

 Proceeding identically as in Definition~\ref{expand},  we define an
$\{\cF_t, \ 0 \leq t \leq T\}$-consistent expectation
$\{\widetilde\e[\cdot|\f_t],\ 0 \leq t \leq T\}$ by extending
$\e_{t,T}^{+}$ from $L^2_+(\cF_T)$ to $L^2(\cF_T)$.

\bl~\label{ee'} Let $\{\e_{s,t}[\cdot], 0\le s\le t\le T\}$ be an
$\f$-consistent dynamic operator with the negative floor $S\le 0$
and satisfy (H1) and (H2). For $X\in \lfst$,

(i) we have for $t\in [0,T]$,
\begin{equation}~\label{equ} \widetilde\e [X|\f_t]\leq
\e_{t,T}[X].
\end{equation}
Furthermore, if $\e_{t,T}[X]>S_t\ \as$ for any $t\in[r,T]$ with
some $r\ge 0$, then
$$\widetilde\e [X|\f_t]=\e_{t,T}[X]\q  a.s. \hbox{ \rm for any }
t\in[r,T]. $$

(ii) The process $\{\cE[X|\cF_t], 0\le t\le T\}$ is an
$\widetilde\cE$-supermartingale. \el

{\bf Proof.} First, we prove assertion (i). Set $X_n=X1_{\{X\geq
-n\}}$. From (H2), we have
$$\lim_{n \rightarrow \infty} \e_{t,T}[X_n]=\e_{t,T}[X]. $$
From the definition of $\widetilde\e [\cdot|\f_t]$ and (H1), we
have
\begin{eqnarray}~\label{7}
\widetilde\e [X|\f_t] &=& \lim_{n \rightarrow \infty}(\e
[X_n+n|\f_t]-n) \nonumber\\  &=& \lim_{n \rightarrow
\infty}(\e_{t,T}[X_n+n]-n)\nonumber\\
 &\leq& \lim_{n \rightarrow
\infty}\e_{t,T}[X_n]=\e_{t,T}[X].
\end{eqnarray}
Therefore, we have $$\widetilde\e [X|\f_t]\leq\e_{t,T}[X].
$$

If $\e_{t,T}[X]>S_t \ \as $ for $t\in[r,T]$, then we have from
(D1) that
$$\e_{t,T}[X_n]\ge \e_{t,T}[X]>S_t \ \as \hbox{ \rm for } t\in[r,T].
$$
In view of (H2), we see that the inequality in (\ref{7}) is
actually an equality for $t\in [r,T]$. Therefore,
$$\widetilde\e [X|\f_t]=\e_{t,T}[X] \ \as \hbox{ \rm for } t\in [r,T].
$$

We now prove assertion (ii). Set $X_t=\e_{t,T}[X]$ for $t\in
[0,T]$. Then from assertion (i), we have for all $0\leq s \leq t
\leq T$
$$
\widetilde\e [X_t|\f_s]\leq \e_{s,t}[X_t]= X_s.
$$
Therefore, $\{\e_{t,T}[X], 0\le t\le T\}$ is an $\widetilde\e
$-supermartingale.
\endpf

The following nonlinear Doob-Meyer's decomposition theorem, due to
Coquet, Hu, Memin, and Peng~\cite[Theorem 6.3, page 20]{CHMP},
will play a key role in the following arguments.

\bl \label{lmdbmy} Let $\cE[\cdot]$ be an $\{\cF_t, 0\le t\le
T\}$-consistent expectation, and let $\{Y_t, 0\le t\le T\}$ be a
continuous $\cE$-supermartingale such that \be E[\sup_{t\in
[0,T]}|Y_t|^2]<\infty. \ee Then there exists an $A\in
\cL^2_\cF(0,T)$ such that $A$ is continuous and increasing with
$A_0=0$, and such that $Y+A$ is an  $\cE$-martingale.\el

Using Lemmas~\ref{ee'} and~\ref{lmdbmy}, we can show the
following.

\bt~\label{expand_consistent}  Let $\{\e_{s,t}[\cdot], 0\le s\le
t\le T\}$ be an $\f$-consistent dynamic operator with the negative
floor $S\le 0$ and satisfy (H1) and (H2). Let
$\e_{t,T}^{+}[\cdot]$ denote the restriction of $\e_{t,T}[\cdot]$
on $L^2_+(\cF_T)$, and $\widetilde\e [\cdot |\f_t]$ be the
extension of $\e_{t,T}^{+}[\cdot]$ given in Section 2. Then for
$X\in L^2(\cF_T;S_T)$, there exists an $\{\f_t, 0\le t\le
T\}$-adapted continuous increasing process $\{A_t^X, 0\le t\le
T\}$ with $A_0^X=0$ such that
\begin{equation}
\e_{t,T}[X]=\widetilde\e [X+A_T^X|\f_t]-A_t^X \hbox{ \it for any }
t\in [0,T].
\end{equation}
\et

{\bf Proof.} Assertion (ii) of Lemma~\ref{ee'} shows that
$\{\cE_{t,T}[X], 0\le t\le T\}$ is an $\widetilde
\cE$-supermartingale.

In version of (D2) and (H1), we have
$$ S_t\le\cE_{t,T}[X]\le \cE_{t,T}[X]-S_t=\cE_{t,T}[X]-\cE_{t,T}[S_t]\le \cE^\mu_{t,T}[X-S_t].$$
Therefore, we have for some positive constant $C_1$,
$$
\ba{rcl} &&\displaystyle E\sup_{t\in [0,T]}|\cE_{t,T}[X]|^2\\
&\le&\displaystyle  E\sup_{t\in
[0,T]}\cE_{t,T}^\mu[X-S_t]^2+E\sup_{t\in
[0,T]}|S_t|^2 \\
&\le&\displaystyle  E\sup_{t\in [0,T]}\cE_{t,T}^\mu
[|X|]^2+2E\sup_{t\in
[0,T]}|S_t|^2 \\
&\le&\displaystyle  C_1 E[X^2]+E\sup_{t\in [0,T]}|S_t|^2 <\infty.
\ea$$

In view of Lemma~\ref{lmdbmy}, there exists a continuous
increasing $\{\f_t, 0\le t\le T\}$-adapted process $\{A_t^X, 0\le
t\le T\}$ with $A_0^X=0$ such that $\{\cE_{t,T}[X]+A_t^X, 0\le
t\le T\}$ is an $\widetilde\e $-martingale. Therefore,
$$\widetilde\e [\cE_{T,T}[X]X+A_T^X|\f_t]
=\e_{t,T}[X]+A_t^X, \q \forall t\in [0,T].$$ While
$\cE_{T,T}[x]=X$, the proof is then complete. \endpf

 In the following, we introduce three lemmas for subsequent arguments.

 Identical to the case of an $\{\cF_t, 0\le t\le
T\}$-consistent expectation, we can show

\bl Let $\{\e_{s,t}[\cdot], 0\le s\le t\le T\}$ be an
$\f$-consistent dynamic operator with the negative floor $S\le 0$
and satisfy (H1) and (H2). For all $ 0\leq s\leq t \leq T$,
$X,Y\in L^{2}(\mathcal{F}_t;S_t)$, and $B\in \f_s$, we have
\begin{equation}
\e_{s,t}[X1_B+Y1_{B^c}]=\e_{s,t}[X]1_B+\e_{s,t}[Y]1_{B^c}.
\end{equation}
\el

\bl~\label{da0} Let $\{\e_{s,t}[\cdot], 0\le s\le t\le T\}$ be an
$\f$-consistent dynamic operator with the negative floor $S\le 0$
and satisfy (H1) and (H2). Let the process $\{A_t^X, 0\le t\le
T\}$ with $A_0^X=0$ be given in Theorem~\ref{expand_consistent}.
For a given $X\in \lfst$ and some $r\in [0,t)$, if
$$\e_{s,T}[X](\omega)>S_s(\omega), \q s\in[r,t], \as \omega\in
B, B\in \f_r,
$$then
\begin{equation}
A_s^X1_{B}=A_t^X1_{B},\q \as \hbox{ \rm for } s\in[r,t].
\end{equation}
\el

 {\bf Proof. } For $\epsilon>0$, set
$Y_t=\e_{t,T}[X]1_B+\epsilon1_{B^c}.$ Then
$$\e_{s,t}[Y_t]>S_s,\ s\in[r,t].$$
From Lemma \ref{ee'}, $Y_t$ satisfies
\begin{equation}
\e_{s,t}[Y_t]=\widetilde\e [Y_t|\f_s].
\end{equation}
Then  from Lemma \ref{da0}, we have for $\forall s\in[r,t]$,
\begin{eqnarray}~\label{11}
\widetilde\e [Y_t+(A_t^X-A_s^X)1_B|\f_s] &=& \widetilde\e
[\widetilde\e
[X+A_T^X-A_t^X|\f_t]+(A_t^X-A_s^X)|\f_s]1_{B}+\epsilon1_{B^c}
\nonumber
\\
&=&\widetilde\e [\widetilde\e
[X+A_T^X-A_s^X|\f_t]|\f_s]1_{B}+\epsilon1_{B^c} \nonumber
\\
&=&\widetilde\e [X+A_T^X-A_s^X|\f_s]1_{B}+\epsilon1_{B^c}\nonumber
\\ &=& \e_{s,T}[X]1_B+\epsilon1_{B^c},\ \nonumber
\\ &=& Y_s,\ \  a.s.
\end{eqnarray}
and
\begin{eqnarray}~\label{12}
\widetilde\e [Y_t|\f_s]&=&\e_{s,t}[Y_t]\nonumber \\
&=& \e_{s,t}[\e_{t,T}[X]]1_B+\epsilon1_{B^c} \nonumber \\
&=& Y_s, \q a.s..
\end{eqnarray}

 Since $A_t^X1_{B}-A_s^X1_{B}\geq 0$ , we have from (\ref{11}) and (\ref{12})
that
$$A_s^X1_{B}=A_t^X1_{B},\q \as \hbox{ \rm for } \forall s\in [r,t]. $$
\endpf

\bl ~\label{min_At} Let $\{\e_{s,t}[\cdot], 0\le s\le t\le T\}$ be
an $\f$-consistent dynamic operator with the negative floor $S\le
0$ and satisfy (H1) and (H2). Let the process $\{A_t^X, 0\le t\le
T\}$ with $A_0^X=0$ be given in Theorem~\ref{expand_consistent}.
For a given $X\in \lfst$, $\{\e_{t,T}[X],0\leq t \leq T\}$
satisfies the following
\begin{equation}
\int_0^T(\e_{t,T}[X]-S_t)dA_t^X=0, \ \as.
\end{equation}
\el

 {\bf Proof.} We firstly prove
$1_{\{\e_{t,T}[X]>S_t\}}dA_t^X=0\ a.s.$. Set
$$B_n:=\{\omega:\e_{r,T}[X](\omega)>S_r(\omega), \forall [(t-1/n)\vee 0]<r<[(t+1/n)\wedge
T])\};$$
$$C_n:=\cup\{C\in \f_{(t-1/n)\vee 0}:C\subseteq B_n\}.$$
Since $\{\mathcal{F}_t,0 \leq t \leq T\}$ is the natural
filtration of ${B_t}$, augmented by all $P$-null sets of
$\mathcal{F}$, we have $$
\cup_{n=1}^{\infty}C_n=\{\e_{t,T}[X]>S_t\}.$$
% \\Or else, there exists a subset $C\in \f_t$ satisfying $$P(C)>0\ \mbox{and}\ C\notin \f_s,\forall s<t$$
Lemma~\ref{da0} implies that $1_{C_n}dA_t^X=0,\as .$ Therefore,
$$1_{\{\e_{t,T}[X]>S_t\}}dA_t^X=0\q  \as .$$
For $\e_{t,T}[X]-S_t\geq 0$, we have \be
\int_0^T(\e_{t,T}[X]-S_t)dA_t^X =
\int_0^T(\e_{t,T}[X]-S_t)1_{\{\e_{t,T}[X]>S_t\}}dA_t^X =0. \ee
\endpf

\bt~\label{find} Theorem~\ref{RBSDEConverse} is true in the case
of the negative floor $S\le 0$. \et

{\bf Proof.} From Theorem~\ref{BSDEConverse}, there exists a
function $g=g(t,z):\Omega \times [0,T]\times \mathbb{R}^d $
satisfying (C2)-(C3) and $g(\cdot,\cdot,0)\equiv 0$, such that the
following holds:
$$\widetilde\e [Y|\cF_t]=\e^g[Y|\cF_t],\ \forall Y\in L^2(\f_T), t\in [0,T].$$
Therefore, for $X\in \lfst$, we have $\widetilde\e
[X+A_T^X|\cF_t]=\e^g[X+A_T^X|\cF_t]$ for $t\in [0,T]$.

From the definition of $\cE^g[X+A_T^X]$ and
$\cE^g[X+A_T^X|\cF_t]$, we know that there is unique $Z\in
\cL^2_\cF(0,T;\mathbb{R}^d)$ such that
$$\e^g[X+A_T^X|\cF_t]=X+A_T^X+\int_t^Tg(s,Z_s)ds-\int_t^T\langle Z_s, dB_s\rangle, \ \as \hbox{ \rm for any } t\in [0,T].$$

 Set
\begin{equation}
\widetilde X_t:=X+\int_t^Tg(s,Z_s)ds-\int_t^TZ_s\,
dB_s+A_T^X-A_t^X.
\end{equation}
Then from Lemma~\ref{expand_consistent}, we have
\begin{eqnarray}
\e_{t,T}[X] &=& \widetilde\e [X+A_T^X|\f_t]-A_t^X\nonumber \\
&=&\cE^g[X+A_T^X|\cF_t]-A_t^X=\widetilde X_t.
\end{eqnarray}
Since $\widetilde X_t=\e_{t,T}[X]\geq S_t$, it follows from
Lemma~\ref{min_At} that $\{(\widetilde X_t,Z_t,A_t^X),0\leq t\leq
T \}$ is the solution of RBSDE $ (X,g,S)$. That is,
\begin{equation}
\e_{t,T}[X]=\widetilde X_t=\mathcal{E}_{t,T}^{r;g,S}[X],\ \
\forall X\in \lfst.
\end{equation}
\endpf

\br~\label{unique} From the proof of Theorem~\ref{find},
$\{(\e_{t,T}[X],Z_t,A_t^X), 0\leq t\leq T \}$ is the solution of
RBSDE $(X,g,S)$. Therefore, the increasing process $\{A_t^X, 0\le
t\le T\}$ is unique for a given $X\in \lfst$. \er

\section{The general case of the upper bounded floor: the proof of Theorem~\ref{RBSDEConverse}.}

Now consider the general case of the upper bounded floor $S$.

\bde Define a new dynamic operator
\begin{equation}
\e_{s,t}^C[\cdot]:L^{2}(\mathcal{F}_t;S_t-C)
 \longrightarrow L^{2}(\mathcal{F}_s;S_s-C), 0\leq s\leq t\leq
 T
\end{equation}
by
\begin{equation}
\e_{s,t}^C[X]:=\e_{s,t}[X+C]-C, \forall X\in
L^{2}(\mathcal{F}_t;S_t-C).
\end{equation}
\ede

It is easy to prove that $\e_{s,t}^C[X]$ for  $ X\in
L^{2}(\mathcal{F}_t;S_t-C),$ satisfies all the conditions (H1),
(H2), and (D1)-(D4). Then, applying Theorem~\ref{find} to
$\e_{s,t}^C[X]$, we have Theorem~\ref{RBSDEConverse}.

\vspace{1cm} \noindent {\bf Acknowledgment. } Both authors thank
Professor Shige Peng for his helpful comments.

\end{document}